\documentclass[10pt, letterpaper]{article}

\usepackage[margin=1in]{geometry} 

\usepackage{amsmath,amssymb,amsfonts, amsthm}
\usepackage[]{algpseudocode}\usepackage{graphicx} 
\usepackage{nicefrac}
\newtheorem{theorem}{Theorem}[section]

\newtheorem{proposition}[theorem]{Proposition}

\newtheorem{remark}[theorem]{Remark}

\newcommand{\lb}{\left\{}
\newcommand{\rb}{\right\}}

\def\thickhline{\noalign{\hrule height1pt}}

\newcommand{\RR}{\mathbb{R}}

\newcommand{\abs}[1]{\lvert #1 \rvert}

\usepackage{amsopn}

\usepackage{algorithm}
\usepackage{textcomp}
\usepackage{xcolor}
\usepackage{graphics}
\usepackage{caption}
\usepackage[hidelinks]{hyperref}
\usepackage{subcaption}
\usepackage{arydshln}
\usepackage{mathtools}
\usepackage{relsize}
\usepackage{BOONDOX-uprscr}
\usepackage{dirtytalk}

\usepackage{color}
\definecolor{darkgreen}{rgb}{0,0.5,0}
\usepackage[colorinlistoftodos]{todonotes}

\usepackage{soul}
\setulcolor{red} 
\setstcolor{red}
\sethlcolor{yellow}

\usepackage{amsmath,amssymb,amsfonts,amsthm}
\usepackage{graphicx}
\usepackage{nicefrac}
\usepackage{multicol}
\usepackage{tikz}
\newcommand\anno{0} 
\newcommand\target{5} 
 
\pgfmathsetmacro\myend{\target-1-\anno}
\pgfmathsetmacro\tix{1}
\pgfmathsetmacro\myspacing{5.7/(\target-1-\anno)}

\title{\LARGE \bf Energy Efficient Nonlinear
Microscopic Dynamical Model for Autonomous and Electric Vehicles}

\author{
Yuneil Yeo$^{1}$, Jaewoong Lee$^{1}$, Scott Moura$^{1}$, and Maria Laura Delle Monache$^{1}$\vspace{0.5em}\\
\small $^{1}$Department of Civil and Environmental Engineering, University of California, Berkeley, CA, USA\\
\small {\tt yuneily@berkeley.edu, ljw7696@berkeley.edu, smoura@berkeley.edu, mldellemonache@berkeley.edu}\\[0.5em]
\small Corresponding author: Yuneil Yeo ({\tt yuneily@berkeley.edu})
}

\date{}

\begin{document}

\maketitle
\thispagestyle{empty}
\pagestyle{empty}

{
\renewcommand\thefootnote{}
\footnotetext{© 2025 IEEE. Personal use of this material is permitted.  Permission from IEEE must be obtained for all other uses, in any current or future media, including reprinting/republishing this material for advertising or promotional purposes, creating new collective works, for resale or redistribution to servers or lists, or reuse of any copyrighted component of this work in other works.}

\addtocounter{footnote}{0}
\renewcommand\thefootnote{\arabic{footnote}}
}

\begin{abstract}
This article proposes a nonlinear microscopic dynamical model for autonomous electric vehicles (A-EVs) that considers battery energy efficiency in the car-following dynamics. The model builds upon the Optimal Velocity Model (OVM), with the control term based on the battery dynamics to enable thermally optimal and energy-efficient driving. We rigorously prove that the proposed model achieves lower energy consumption compared to the Optimal Velocity Follow-the-Leader (OVFL) model. Through numerical simulations, we validate the analytical results on the energy efficiency. We additionally investigate the stability properties of the proposed model. 
\end{abstract}

\section{Introduction}\label{sec:introduction}

The recent emergence of autonomous electric vehicles (A-EVs) is expected to enhance traffic flow and safety with V2X communications. However, the development of A-EVs introduced new challenges in designing their dynamic systems. While the traditional microscopic models focused on traffic safety and flow dynamics, A-EV models must additionally consider energy consumption and thermal management to effectively maintain battery health, extend driving range, and improve performance reliability \cite{bandhauer2011critical, hwang2024review}. These new requirements necessitate the need to incorporate energy efficiency into the modeling and control of vehicle-level behavior. 

Microscopic traffic flow models are well-suited for this purpose, as they focus on traffic patterns at the most granular level and capture the vehicle-level dynamics. Notable microscopic models include Pipes' car following model \cite{pipes1953operational}, Newell’s car-following model \cite{newell1961nonlinear}, and the Optimal Velocity Model (OVM) \cite{bando1995dynamical}. Among these, OVM effectively captures complex real-world traffic patterns, such as stop-and-go waves and traffic breakdown, due to its dynamic structure \cite{dong2009velocity}. Building on these advantages, several studies have modified OVM to better reflect driving patterns, leading to extensions like the Optimal Velocity Follow-the-Leader (OVFL) model \cite{10384086}. Some works have developed nonlinear dynamical models for autonomous vehicles to provide flexible control for desired speeds \cite{matin2024analytical} while also enhancing safety and ride comfort \cite{van2006impact, milanes2014modeling, shladover2012impacts}. Investigations into the properties of OVM-based models have examined factors such as linear local stability \cite{ward2011criteria, wilson2011car}, string stability \cite{bando1998analysis, wilson2011car, gunter2020commercially, giammarino2020traffic}, convergence rates \cite{matin2020nonlinear, 9147244}, and safety \cite{nick2022near, 10384086}.

Building on these past works, several recent studies have proposed energy-efficient control strategies for microscopic models, particularly focused on A-EVs in mixed traffic environments. For instance, the Energy-Efficient Electric Driving Model (E3DM) has been developed as a rule-based adaptive cruise control scheme for A-EVs, effectively enhancing regenerative braking efficiency by adjusting inter-vehicle spacing \cite{lu2019energy}. Some works consider an OVM-based longitudinal controller design where either the control term considering the motor dynamics is integrated \cite{li2019electric} or where the parameters are determined by solving optimization problems aimed at maximizing energy efficiency \cite{shen2023energy}. Moreover, other researchers have formulated dedicated optimization problems aimed at maximizing energy savings in vehicle operations while considering other aspects like riding comforts, vehicle dynamics in terms of acceleration, and energy consumption behaviors \cite{xian2022economic, chen2020battery, zhang2021energy}. Many works incorporate a Model Predictive Control (MPC) scheme for energy-efficient electric vehicles with the consideration of vehicle dynamics and battery dynamics, including current and state of charge (SOC) \cite{han2020battery, yeom2022model, xie2019predictive}.

However, despite these significant contributions, there remains a notable gap in the literature: an Optimal Velocity Model (OVM)-based nonlinear microscopic dynamical model that integrates control term(s) that consider the battery energy efficiency in terms of cell temperature, heat generation, and power loss through detailed battery dynamics. Addressing this gap is critical for advancing control strategies in A-EVs that balance the driving performance with energy efficiency. Since energy consumption directly affects driving range, battery health, and thermal performance, incorporating these factors into vehicle-level models is key to ensuring reliable, efficient, and sustainable operation. The development of such models, therefore, ultimately supports the deployment of energy-aware autonomous electric vehicles. 

\textbf{Focus and Contribution.} In this paper, we propose a microscopic traffic model for A-EVs that achieves energy-efficient driving. Our approach builds on an OVM-based framework by integrating comprehensive battery dynamics—including state-of-charge, thermal behavior, and internal resistance—directly into the traffic model. This integration allows for real-time adjustments of energy parameters such as battery temperature and heat generation, thereby ensuring that the system operates under optimal thermal conditions while maximizing energy efficiency.

The paper is organized as follows. In Section \ref{sec: model}, we present the nonlinear microscopic dynamical model for autonomous electric vehicles (A-EVs), with a focus on battery energy efficiency. Section \ref{sec: battery_efficiency} details the battery dynamics, explores their relationship with vehicular dynamics, and explains the formulation of the control terms in the model. In Section \ref{sec: performance}, we present the analytical result that the energy consumption of a vehicle under a proposed model is lower than the energy consumption of a vehicle under OVFL. In Section \ref{sec: numerical_simulation}, we present numerical simulations to validate the analytical work on the improved energy efficiency of the proposed model over OVFL. Finally, Section \ref{sec:stability} discusses the stability of the proposed model. 

\section{Mathematical Model}\label{sec: model} 
We first propose the nonlinear microscopic dynamical model in the following form:
\begin{equation}\label{E:main_dynamics}
\begin{cases}
\dot x_l  = v_l \\
\dot v_l = a_l \\
\dot x = v\\
\dot v = \alpha \lb V(x_l-x)-v \rb + \beta\frac{v_l-v}{(x_l-x)^2} - \kappa v^2 \frac{(v_l-v)^2}{(v_l-v)^2+\epsilon} \\
(x_l(t_\circ) , x(t_\circ) , v_l(t_\circ), v(t_\circ)) = (x_{l, \circ}, x_\circ, v_{l, \circ}, v_\circ), 
\end{cases}
\end{equation}

where \(t_f > t_\circ\) is a fixed time horizon, \(t_\circ\) is the initial time, and \(\epsilon\) is a small positive number (e.g., \(10^{-6}\)) introduced to preserve the smoothness of the model. The leading vehicle’s position and velocity are denoted by \(x_l\) and \(v_l\), respectively, while the following (ego) vehicle’s position and velocity are denoted by \(x\) and \(v\). The leading vehicle's acceleration, \(a_l\), should be chosen such that its velocity remains non-negative and does not exceed the maximum allowable velocity $v_{\max}$.

The optimal velocity function \(V\) is smooth, monotonically increasing, bounded, and Lipschitz continuous. For this study, we consider
\begin{equation}  \label{E:ovm_function}
{V}(u) = \tanh(u-2) + \tanh(2), \quad \text{for any } u \in \RR.
\end{equation}
For properties of this function, see \cite{nick2022near, bando1995dynamical}, and for a sensitivity analysis with respect to the non-negative coefficients \(\alpha\) and \(\beta\), refer to \cite{10384086}.  The non-negative coefficient \(\kappa\) represents the weight on the energy efficiency of autonomous electric vehicles (A-EVs), capturing the trade-off between achieving higher velocity and minimizing energy loss. The parameter values should be chosen such that $\beta > \alpha$ to avoid collision and that  $\kappa$ is smaller than both $\alpha$ and $\beta$ to maintain stability.

\begin{figure}[!htb]  \vspace{0.2cm}
  \begin{subfigure}[b]{0.5\linewidth}
    \centering
    \includegraphics[width=0.9\linewidth]{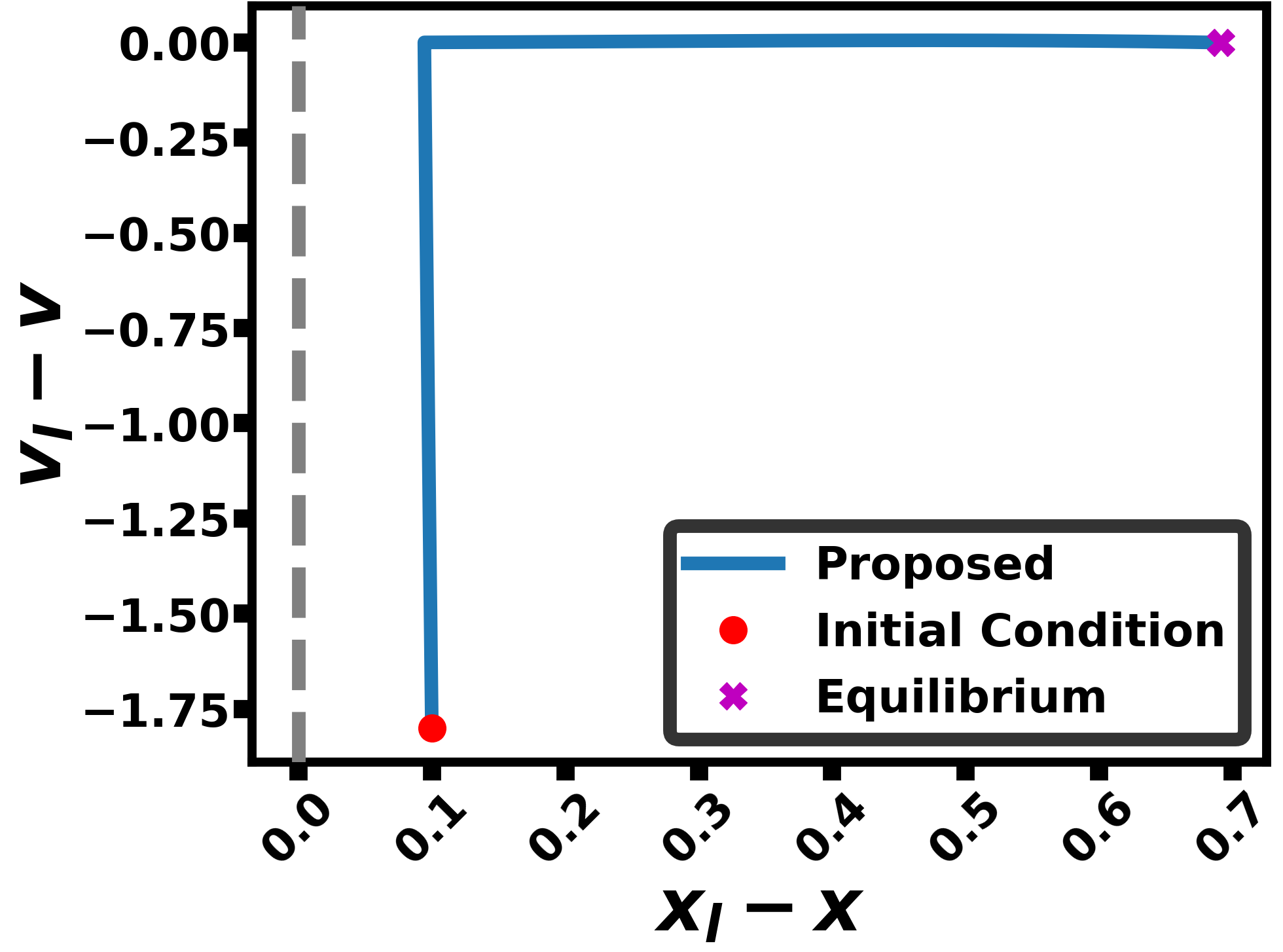} 
    \caption{$x_l$ = 0.1, $v$ = 1.9, $v_l$ = 0.1} 
    \label{fig: left_top_constant_velocity} 
  \end{subfigure}
  \begin{subfigure}[b]{0.5\linewidth}
    \centering
    \includegraphics[width=0.9\linewidth]{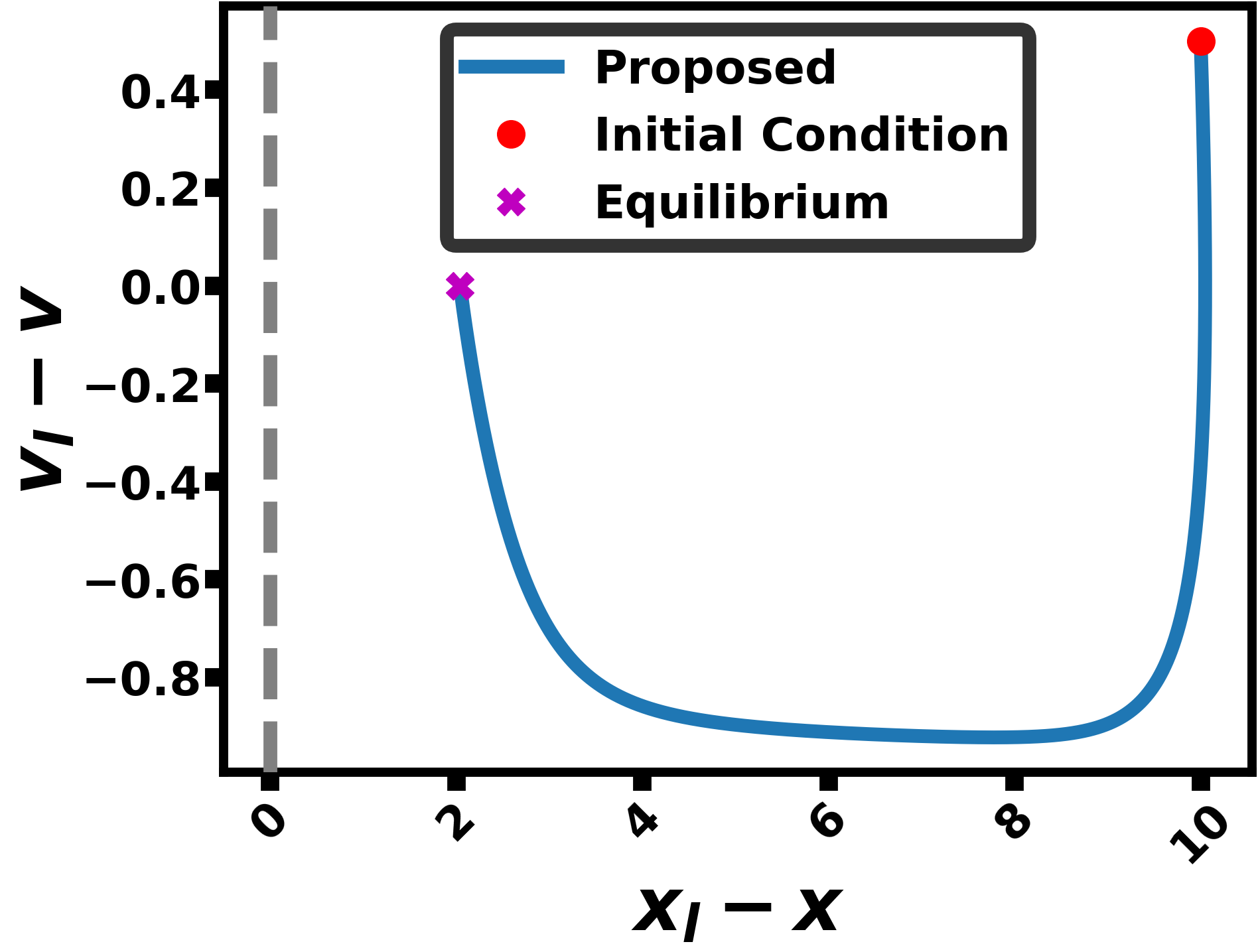} 
    \caption{$x_l$ = 10, $v$ = 0.5, $v_l$ = 1} 
    \label{fig: right_top_constant_velocity} 
  \end{subfigure} 
  \caption{\small{Trajectory of the proposed dynamics. The initial condition is shown as a red dot, while the equilibrium is shown as a magenta x mark.  Parameters: $\alpha$ = 2, $\beta$ = 3, $\kappa$ = 0.03, $x$ = 0, Simulation Time = 700, $a_l$ = 0}}
  \label{fig: traj_proposed} 
\end{figure}

Figure \ref{fig: traj_proposed} shows the trajectory of the proposed model under different initial conditions. These plots illustrate how the spacing and relative velocity between the leading and ego vehicles evolve over time under the proposed model. 

To further examine the proposed model, we compare the proposed model with the existing microscopic traffic model for AVs: OVFL. While $\dot x_l, \dot v_l, \dot x$ of OVFL are the same as the proposed model, $\dot v$ of OVFL is \begin{equation*}
\dot{v} = \alpha \left[ V(x_l - x) - v \right] + \beta \frac{(v_l - v)}{(x_l - x)^2}. 
\end{equation*}

We now plot OVFL and the proposed model to check if the proposed model yields reasonable behaviors compared to OVFL.
\begin{figure}[!htb] 
  \begin{subfigure}[b]{0.5\linewidth}
    \centering
    \includegraphics[width=0.9\linewidth]{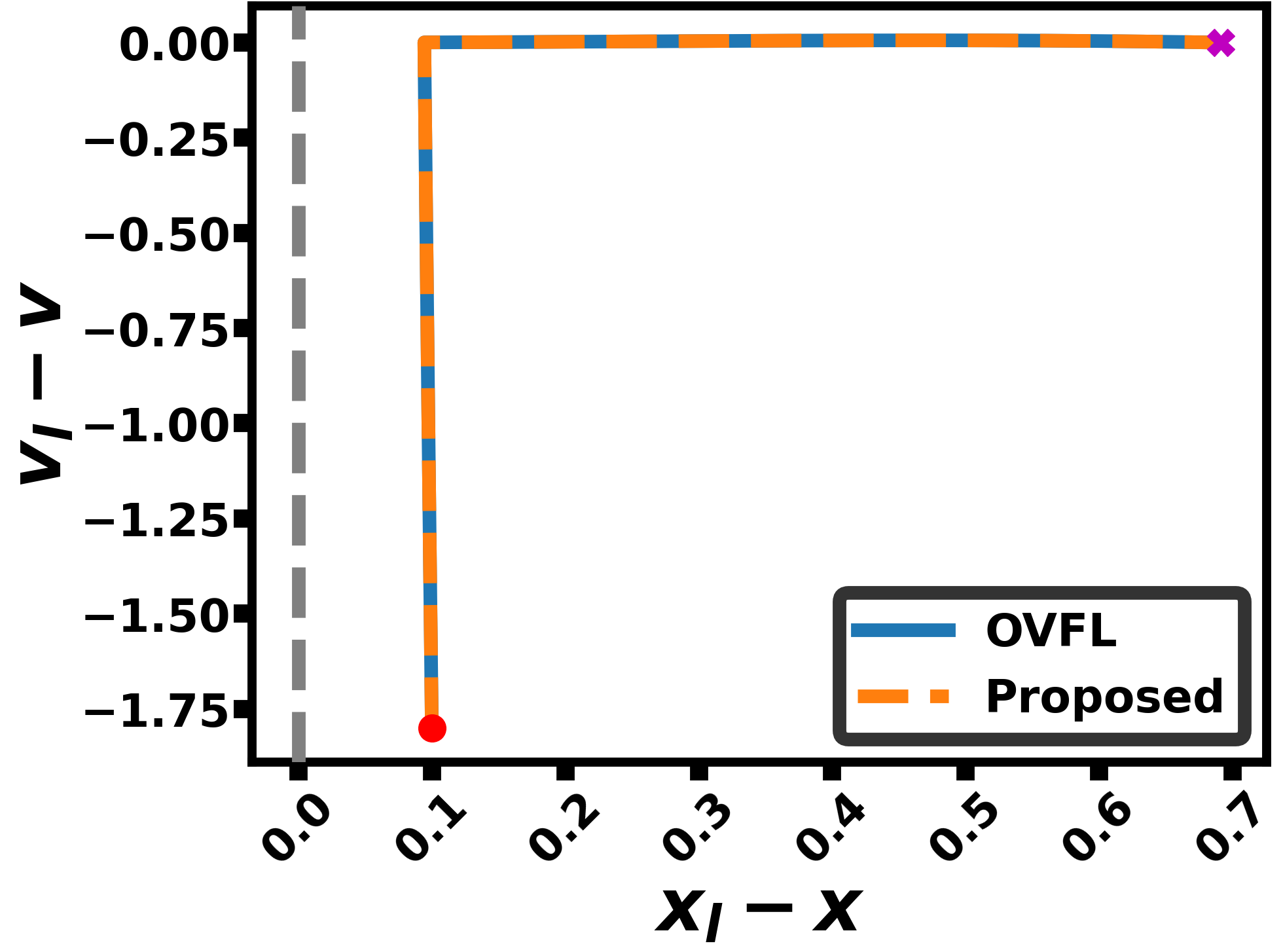} 
    \caption{$x_l$ = 0.1, $v$ = 1.9, $v_l$ = 0.1}
    \label{fig: left_top_constant_velocity} 
  \end{subfigure}
  \begin{subfigure}[b]{0.5\linewidth}
    \centering
    \includegraphics[width=0.9\linewidth]{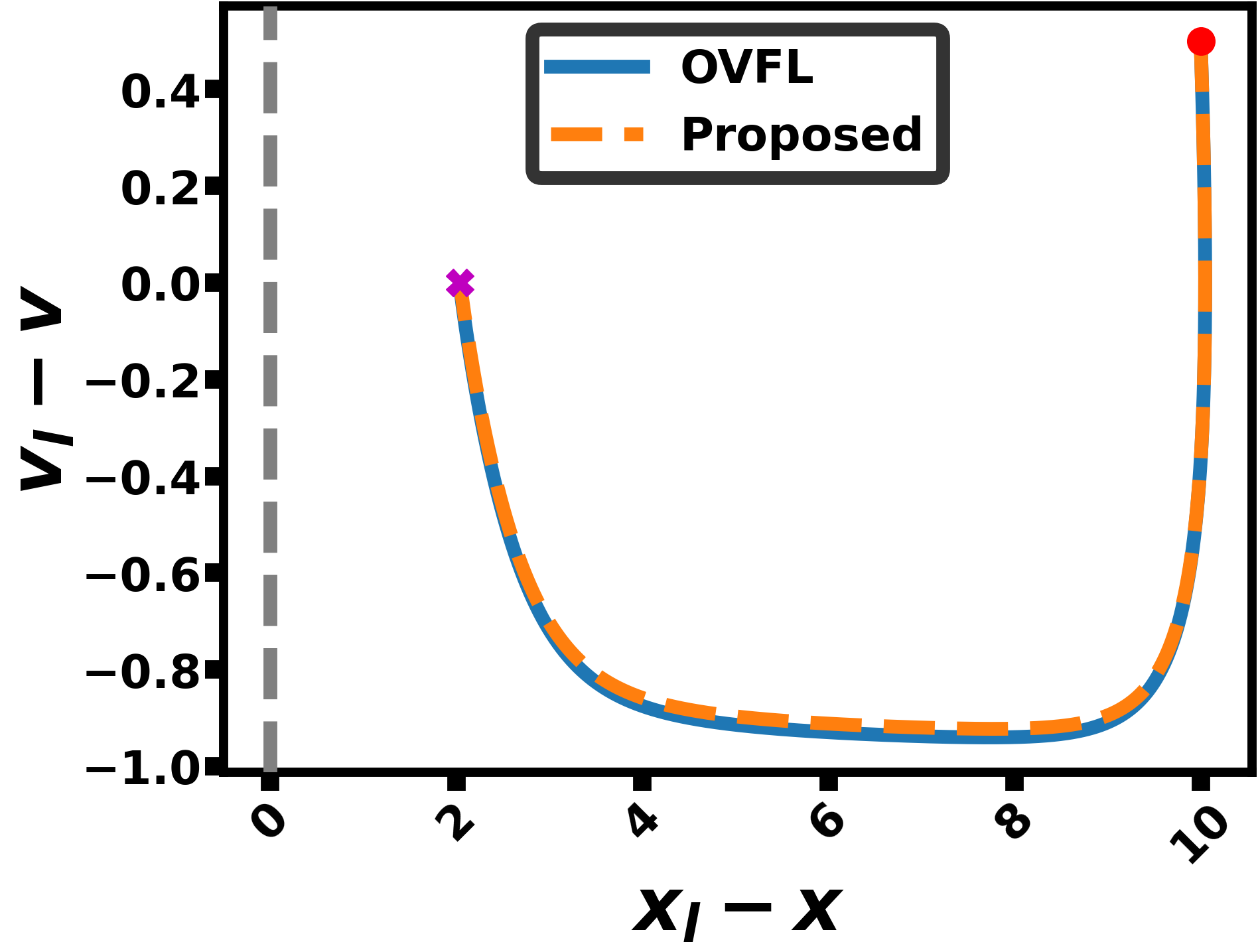} 
    \caption{$x_l$ = 10, $v$ = 0.5, $v_l$ = 1} 
    \label{fig: right_top_constant_velocity} 
  \end{subfigure}
  \caption{\small{Comparison between the trajectory of OVFL and the proposed dynamics. The initial condition is shown as a red dot, while the equilibrium is shown as a magenta x mark. Parameters: $\alpha$ =2, $\beta$ = 3, $\kappa$ = 0.03, $x$ = 0, Simulation Time = 700, $a_l$ = 0}}
  \label{fig:comparison} 
\end{figure}

Figure \ref{fig:comparison} shows that the trajectories of the proposed model closely match those of OVFL. Specifically, both systems reach equilibrium without crossing the gray dashed line, which serves as an indicator of a collision. This behavior suggests that the proposed model retains collision-avoidance properties of OVFL while offering improved energy efficiency through the additional control term. The battery efficiency of the proposed model is discussed in a later section.

The properties of the model in the case of two consecutive vehicles can be extended to the platoon case as the leading vehicle's acceleration $a_l$ is modeled in a general form and is not tied to a specific model. This independence allows the two-vehicle interaction structure to be recursively applied throughout the platoon \cite{nick2022near}. 

\section{Battery Efficiency} \label{sec: battery_efficiency}

The primary motivation of this paper is to develop a dynamic nonlinear model for autonomous electric vehicles that explicitly integrates battery dynamics and examines the interplay between battery and vehicular dynamics to account for battery efficiency. To this end, we begin by outlining several key assumptions underlying the battery dynamics used in our model formulation:

\begin{enumerate}
\item \textbf{Short-duration trip}
\item \textbf{Constant battery temperature:} In a short driving duration, it is reasonable to assume a steady inlet air coolant, resulting in constant air temperature and flow, which makes the battery temperature time-invariant.
\item \textbf{Constant internal resistance:} With a constant battery temperature, internal resistance is treated as time-invariant.
\item \textbf{Single-cell representation:} 
Battery packs are usually composed of identical battery cells. Therefore, the analysis on the single cell simplifies the analysis without loss of generality. 
\item \textbf{Linear motor-battery power conversion relationship:} A linear relationship simplifies modeling while capturing essential behavior.
\item \textbf{Exclusion of entropic heat generation:} In a short-duration trip, entropic heat is negligible compared to resistive heating.
\item \textbf{Neglecting battery aging and degradation:} Battery aging and degradation are minimal in a short-duration trip. Therefore, these effects can be excluded.
\item \textbf{Fixed and known battery and coolant properties:} For a short-duration trip, these properties are typically well-defined, supporting a precise analysis.
\end{enumerate}

Based on these assumptions, we now consider the battery dynamics using a single-cell circuit model. From \cite{perez2012parameterization}, the battery terminal voltage $V_T$ is defined as 
\begin{equation} \label{E:terminal_voltage}
    V_T = V_{OCV} - V_1 - V_2 - IR_s
\end{equation}
where $V_{OCV}$ is the open circuit voltage, and $I$ is the current (the input to the battery cell). The current $I$ is defined as positive during discharge and negative during charging, including scenarios involving regenerative braking. $R_s$ is the internal resistance of the cell, where it is assumed to be a time-invariant value based on assumptions.

$V_1$ and $V_2$ evolves according to the following dynamics: 
\begin{equation}\label{E:V_1_V_2}
    \begin{cases}
        \dot V_1 = -\frac{1}{R_1C_1}V_1 + \frac{1}{C_1}I,\\
        \dot V_2 = -\frac{1}{R_2C_2}V_2 + \frac{1}{C_2}I,\\
    \end{cases}
\end{equation}
where formulas of $R_1, R_2, C_1, C_2$ are described in \cite{perez2012parameterization}.

The value of $V_{OCV}$ depends on the state of charge. The state of the charge, $S \in [0,1]$, evolves according to
\begin{equation} \label{E:change_SOC}
    \dot{S} = -\frac{1}{3600C_n} I
\end{equation}
where $C_n$ is the nominal capacity (a constant) in units of Ah. The nominal capacity refers to the amount of electric charge that a new battery can store or provide at 25 °C. 

The cell temperature changes due to heat generation from power loss $Q$, which is 
\begin{equation} \label{E:power_lossQ}
    Q = I(V_{OCV}-V_T).
\end{equation}
The internal power loss $Q$ remains non-negative due to resistive and electrochemical effects, regardless of charging or discharging. 

For energy efficiency, we aim to minimize the rate of power loss $Q$. To investigate how to minimize the rate of power loss $Q$, we set up the finite-time optimal control problem over a finite time interval $t \in [t_\circ, t_f]$:
\begin{equation} \label{E:opt_prob}
\begin{aligned}
    &\min_{I}  \int_{t_\circ}^{t_f} Q(t) dt  = \min_{I}  \int_{t_\circ}^{t_f} I(t) (V_{OCV}(t) -V_T(t)) dt \\ 
    &\text{subject to } \eqref{E:terminal_voltage}, \eqref{E:V_1_V_2}, \eqref{E:change_SOC}, \eqref{E:power_lossQ}, \\ 
        & \quad Q(t) \geq 0, \quad S(t) \in [0,1],\\
        & \quad \text{Constants } R_1, R_2, C_1, C_2, C_n, R_s > 0,\\
        & \quad V_{OCV}(t_\circ) = V_{OCV, \circ}, \quad V_{T}(t_\circ) = V_{T, \circ}, \\ 
        & \quad V_1(t_\circ) = V_{1,\circ}, \quad V_2(t_\circ) = V_{2, \circ}, \quad S(t_\circ) = S_\circ \\ 
\end{aligned}
\end{equation}
where $I(t)$ acts as the control input. 

\begin{proposition}[\textbf{Optimal Current for Minimizing Power Loss}]\label{T:minimizeI}
The minimizer of the optimization problem in \eqref{E:opt_prob} is 
\begin{equation}
    I^*(t) = 0, \quad \forall t \in  [t_\circ, t_f]
\end{equation}
\end{proposition}
\begin{proof}
    We start by expanding \eqref{E:power_lossQ} by plugging \eqref{E:terminal_voltage} in. 
    \begin{equation}
        Q(t) = I(t)(V_1(t) + V_2(t) + R_s \ I(t)).
    \end{equation}
    Then, the objective function of the optimization problem becomes 
    \begin{equation}\label{E:objective_func}
        \min_{I}  \int_{t_\circ}^{t_f} I(t)(V_1(t) + V_2(t)) + R_s \ I^2(t) dt
    \end{equation}

    Since \eqref{E:V_1_V_2} are linear-time invariant Ordinary Differential Equations with $I(t)$ as the input, $V_1(t)+ V_2(t)$ is a linear functional of $I(t)$. Therefore, $I(t)(V_1(t) + V_2(t))$ is a bilinear term while $ R_s \ I^2(t)$ is a quadratic term in terms of $I(t)$. Together, the objective function becomes a quadratic functional of $I(t)$, and its value increases with the magnitude of $I(t)$. So, the minimum can be achieved when the magnitude of the current $I(t)$ is minimized.

    We first evaluate the power loss $Q$ over the time interval $t \in [t_\circ, t_f]$ under the condition $I(t) = 0$. Suppose $I(t) = 0$; then, based on $\eqref{E:V_1_V_2}$, \( V_1(t), V_2(t) \to 0 \) exponentially. Since $I(t) = 0$, the power loss $Q$ would be 0 for $t \in [t_\circ, t_f]$ based on \eqref{E:power_lossQ}. Therefore, the integrand in \eqref{E:objective_func} becomes 0. 
    
    Suppose there exists a subinterval $\mathcal T \subseteq [t_\circ, t_f]$ where $I(t) \neq 0$ for $t \in \mathcal T$. On $\mathcal T$, the term $R_s I(t)^2 > 0$. While the term $I(t)(V_1(t) + V_2(t))$ can be negative, the total $Q(t) \geq 0$ due to the physical constraint, being strictly positive on $\mathcal T$. Therefore, when $I(t) \neq 0$ for $t \in \mathcal T$, 
    \begin{equation*}
    \int_{t_\circ}^{t_f} Q(t)\, dt \geq 0 .
    \end{equation*}
    
Hence, $I(t) = 0$ leads to the minimum value of the objective function in \eqref{E:objective_func}, and is therefore a minimizer of the optimization problem. This completes the proof. \end{proof}

\begin{remark}[\textbf{Battery Discharge Rate Benefit}]
Minimizing the magnitude of $I(t)$, ideally $I(t) = 0$, not only minimizes the power loss $Q$, but also minimizes battery drainage rate in a finite time interval. According to \eqref{E:change_SOC}, the change in state-of-charge and the magnitude of the $I(t)$ are in a negative linear relationship. Therefore, minimizing the absolute value of $I$ slows the battery's discharge rate.
\end{remark}

Motivated by these benefits, we now investigate the condition that minimizes the absolute value of the current $I$. To do so, we first look at the total power output of the battery cell, which is
\begin{equation} \label{E:Power_output_current}
    P_{output} = I V_T .
\end{equation}

\subsection{Cell Power to Battery Power}

We now scale up the cell power to represent the battery power. 
\begin{equation}
    P_{battery, output} = (N_s V_T)(N_p I) = N_s N_p P_{output}
\end{equation}
where $N_s$ is the number of cells connected in series and $N_p$ represents the number of parallel strings. 

By solving for $P_{output}$, we get

\begin{equation}
     P_{output} = \frac{P_{battery, output}}{N_s N_p}.
\end{equation}

\subsection{Battery Power to Motor Power}

An electric vehicle converts battery power to motor power. Although the relationship between the battery power to motor power is nonlinear, the magnitude of the nonlinearity is small. Therefore, we can assume that the battery power to motor power conversion is in a linear relationship, where the efficiency of the electric motor $\eta$ (usually 80-90\%) determines the conversion. Therefore,
\begin{equation}\label{E:P_motor}
P_{motor} = 
    \begin{cases}
        \eta P_{battery, output}, \quad \text{if } P_{battery, output} > 0\\
        \frac{P_{battery, output}}{\eta}, \quad \text{if } P_{battery, output} \leq 0\\
    \end{cases}
\end{equation}
where $P_{motor}$ is the mechanical motor power. The case of $P_{battery, output} > 0$ is the case of discharging ($I > 0$) while $P_{battery, output} \leq 0$ is the case of charging ($I < 0)$. In an electric vehicle, power regeneration is possible during braking or deceleration. Based on the relationship between $P_{battery, output}$ and $P_{output}$, we get
\begin{equation} \label{E:p_output_dynamic_motoer1}
P_{motor} = 
    \begin{cases}
        \frac{\eta P_{output}}{N_s N_p}, \quad \text{if } P_{output} > 0\\
        \frac{P_{output}}{\eta N_s N_p}, \quad \text{if } P_{output} \leq 0\\
    \end{cases}
\end{equation}

\subsection{Motor Power to Car Dynamic}

The relationship between the motor power and the car dynamics, like the acceleration and the velocity, is shown in the equation below
\begin{equation}\label{E:P_output_general}
    F_{accel} v = P_{motor} - F_{total} v
\end{equation}
where $v$ is velocity and $F_{accel}$ is the force due to acceleration or deceleration: 
\begin{equation} \label{E:force_accel}
F_{accel} = m a
\end{equation}
where $m$ is the mass of the vehicle.

$F_{total}$ is the total resistive forces like aerodynamic drag (the resistance exerted by air), rolling resistance (force that resists the motion of the tires), and gravitational forces due to the incline. 

Aerodynamic drag is 
\begin{equation} \label{E:air_drag}
F_d = 0.5 \rho A C_d v^2
\end{equation}
where $\rho$ is an air density, $A$ is the frontal area of the vehicle, and $C_d$ is a drag coefficient. 

Rolling resistance is 
\begin{equation}\label{E:rolling}
    F_{rolling} = C_r m \mathtt g
\end{equation}
where $C_r$ is the rolling resistance coefficient (depends on tire material, tread, and inflation) and $\mathtt g$ is the gravitational acceleration (9.81 $m/s^2$). 

The gravitational force due to the incline is
\begin{equation}\label{E:graviational_force}
F_g = m \, \mathtt g \, \sin(\theta)
\end{equation}
where $\theta$ is the angle of the incline.

Considering \eqref{E:force_accel}–\eqref{E:graviational_force}, \eqref{E:P_output_general} can be written as: \begin{equation} \label{E:p_output_dynamic_motoer}
 ma v = P_{motor} - \left(0.5\rho A C_d\, v^2 + C_r m \mathtt g + m \mathtt g \sin\theta\right) v .
\end{equation}

Based on \eqref{E:Power_output_current}, \eqref{E:p_output_dynamic_motoer}, and \eqref{E:p_output_dynamic_motoer1}, Eq.~\eqref{E:p_output_dynamic_motoer} can be rewritten as: 
\begin{equation} \label{E:p_output_dynamic}
 ma v  = 
\begin{cases}
    \begin{split}
    \frac{\eta I V_T}{N_s N_p} - \left(\frac{1}{2}\rho A C_d\, v^2 + C_r m \mathtt g + m \mathtt g \sin\theta\right) v,\\ \quad \text{if } I > 0 
    \end{split}\\
    \begin{split}
    \frac{I V_T}{\eta N_s N_p} - \left(\frac{1}{2}\rho A C_d\, v^2 + C_r m \mathtt g + m \mathtt g \sin\theta\right) v,\\ \quad \text{if } I \leq 0
    \end{split}
\end{cases}
\end{equation}

From \eqref{E:p_output_dynamic}, the relationship between the current \(I\) and the vehicle dynamics is given by
\begin{equation} \label{E:current_dynamics}
\begin{split}
I = \begin{cases}
    \begin{split}
        \frac{\Biggl[\left(\frac{1}{2}\rho A C_d\, v^2 + C_r m \mathtt g + m \mathtt g \sin\theta\right)v + ma v \Biggr]}{\eta N_s N_p V_T},\\ \quad \text{if } I > 0
    \end{split}\\
    \begin{split}
        \frac{\eta \Biggl[\left(\frac{1}{2}\rho A C_d\, v^2 + C_r m \mathtt g + m \mathtt g \sin\theta\right)v + ma v \Biggr]}{N_s N_p V_T},\\ \quad \text{if } I \leq 0
    \end{split}\\
\end{cases}
\end{split}
\end{equation}

Here, the terminal voltage \(V_T\) continuously varies as the current $I$ changes according to \eqref{E:terminal_voltage}. To minimize the absolute value of the expression in \eqref{E:current_dynamics}, it is sufficient to minimize the numerator. Doing so also minimizes the discrepancy between $V_T$ and $V_{OCV}$, and ultimately $Q$. 

To achieve this, we identify the following design requirements for the deceleration control term:
\begin{enumerate}
    \item \textbf{Proportional to $v^2$:} The deceleration should be applied to counteract the resistive forces, such as aerodynamic drag, rolling resistance, and road grade effects, that grow with $v^2$. 
    \item \textbf{Zero deceleration at rest:} When the vehicle is stationary, no deceleration should be applied to prevent undesired behavior such as backward motion.
    \item \textbf{Zero deceleration at equilibrium:} For energy efficiency, no deceleration should be applied when the ego vehicle’s speed matches that of the leading vehicle (\(v = v_l\)), i.e., at steady-state following.
\end{enumerate}

Considering all these requirements, we propose the following control term for energy efficiency: 
\begin{equation}\label{E: control_term}
    -\kappa v^2 \frac{(v_l-v)^2}{(v_l-v)^2 + \epsilon}
\end{equation}
where $\kappa$ determines the strength of the deceleration term, which is proportional to $v^2$, for energy efficiency. The criteria for selecting the value of $\kappa$ are discussed in a later section. This control term is incorporated into the OVFL model to improve the energy efficiency of A-EVs, resulting in the dynamics given by \eqref{E:main_dynamics}.

\section{Performance Evaluation}\label{sec: performance}

To evaluate the energy efficiency of the proposed model compared to other nonlinear microscopic dynamical models, we consider the energy consumption per unit mass over a particular time interval $t \in [t_\circ, t_f]$ \cite{shen2023energy}, defined as follows: 
\begin{equation} \label{E:energy_consumption}
    \omega = \int_{t_\circ}^{t_f} vg(\dot v) dt
\end{equation}
where $g(u)$ is a function that changes based on the engine and powertrain types. In \cite{shen2023energy}, $g(u) = \max\{u, 0\}$ is a rectified linear function (ReLU) for vehicles with internal combustion engines. For electric vehicles, a different $g(u)$ can be considered as studied in \cite{sciarretta2020energy}. 

For autonomous electric vehicles, we set $g(u)$ to be a leaky rectified linear function (LeakyReLU) to consider the energy generation due to the regenerative braking of electric vehicles. Therefore, $g(u)$ is set as follows:

\begin{equation} \label{E:g_function}
    g(u) = \begin{cases}
    \frac{1}{\eta} u, & \text{if } u \ge 0, \\
    \eta u, & \text{if } u < 0,
    \end{cases}
\end{equation}
where we considered $\eta$ for the inefficiencies during regenerative braking and motoring as we did in \eqref{E:P_motor}. 

The energy consumption of the electric vehicle explicitly depends on motor power $P_{motor}$. Specifically, regenerative braking occurs when $P_{motor} < 0$, not simply whenever a vehicle decelerates. For instance, deceleration caused by frictional or gravitational forces does not produce regenerative power. However, for simplicity and analytical tractability, this paper neglects the external forces, thereby assuming that all decelerations directly result in regenerative braking. Although this assumption excludes cases where the deceleration occurs without energy regeneration, it still provides a clear and practical framework for evaluating and comparing the energy efficiency across different models.

In this context, the following theorem establishes a rigorous comparison between two candidate dynamical models for AVs — the proposed model and the OVFL. The result shows that the inclusion of the extra control term in the proposed model leads to a lower energy consumption per unit mass over a finite time interval, thereby highlighting the energy efficiency benefits of the proposed model.

\begin{theorem}[\textbf{Energy Efficiency}]\label{T:main} 
Let $t_f > t_\circ$ and $v_{\max}>0$. Define the state-space 
\[
\mathcal S = \Bigl\{(x_l,x,v_l,v)\in\RR^2\times\RR_+^2:\; v_l,v\in [0,v_{\max}]\Bigr\}.
\]
Consider the dynamical model
\begin{equation}\label{E:dynamics}
\begin{cases}
\dot x_l  = v_l,\\
\dot v_l = a_l,\\
\dot x = v,\\
\dot v = f(x_l,x,v_l,v),
\end{cases}
\end{equation}
where its initial condition is $(x_{l,\circ},x_\circ,v_{l,\circ},v_\circ)\in \mathcal S$ and its solution is $(x_l,x,v_l,v) \in C^1([t_\circ,t_f]; \mathcal S)$. Note that the leading vehicle's dynamics are the same regardless of the choice of function $f$.

We consider two cases for the function $f$:
\[
f_1 = \alpha\Bigl[V(x_l-x)-v\Bigr] + \beta\,\frac{v_l-v}{(x_l-x)^2} - \kappa\,v^2\,\frac{(v_l-v)^2}{(v_l-v)^2+\epsilon},
\]
and 
\[
f_2 = \alpha\Bigl[V(x_l-x)-v\Bigr] + \beta\,\frac{v_l-v}{(x_l-x)^2},
\]
where $V(u)$ defined as \eqref{E:ovm_function}. 

When energy efficiency is measured by the energy consumption per unit mass over a finite time interval $t \in [t_\circ, t_f]$ defined in \eqref{E:energy_consumption}, then 
\[
\omega_{f_1} \leq \omega_{f_2}.
\]
That is, the energy consumption per unit mass over the finite time interval corresponding to $f_1$ is lower than that corresponding to $f_2$. 
\end{theorem}
\begin{proof} As the term $1/(x_l-x)^2$ introduces the singularity as $x_l-x = 0$, $f_1$ and $f_2$ are not globally Lipschitz on $\mathcal S$. We, therefore, perform the regularization by adding $\delta>0$ to the denominator of the singularity term, resulting in the regularized functions \(f_{1,\delta}\) and \(f_{2,\delta}\):
\[
f_{1,\delta}=\alpha\Bigl[V(x_l - x) - v\Bigr] + \beta\frac{v_l - v}{(x_l - x)^2 + \delta} - \kappa\frac{v^2(v_l - v)^2}{(v_l - v)^2 + \epsilon},
\]
\[
f_{2,\delta} = \alpha\Bigl[V(x_l - x) - v\Bigr] + \beta\,\frac{v_l - v}{(x_l - x)^2 + \delta}.\]

Since \((x_l - x)^2 + \delta\) is bounded away from zero by \(\delta\), both \(f_{1,\delta}\) and \(f_{2,\delta}\) are globally Lipschitz on the state space \(\mathcal{S}\). Let $(x_{l,\delta}^{(i)}, x_\delta^{(i)}, v_{l,\delta}^{(i)}, v_\delta^{(i)}) \in C^1([t_\circ,t_f]; \mathcal S)$, for $i=1,2$, be the unique solutions to the regularized functions. 

Note that, at the same state $(x_l,x,v_l,v) \in \mathcal S$, we have
\[
f_{1,\delta}(x_l,x,v_l,v) = f_{2,\delta}(x_l,x,v_l,v) - \kappa\,v^2\,\frac{(v_l-v)^2}{(v_l-v)^2+\epsilon}.
\]

Since the additional term in $f_{1, \delta}$ is less than or equal to 0, for any given state $(x_l,x,v_l,v) \in \mathcal S$, 
\begin{equation} \label{E:compare_at_same_state}
f_{1,\delta}(x_l,x,v_l,v) \le f_{2,\delta}(x_l,x,v_l,v).
\end{equation}

We now define the differences between the two solutions as follows:
\begin{equation*}
    \begin{cases}
        \Delta x_\delta = x_\delta^{(2)}-x_\delta^{(1)}\\
        \Delta v_\delta = v_\delta^{(2)}-v_\delta^{(1)} = \Delta \dot x_\delta
    \end{cases}
\end{equation*}

In addition, 
\begin{equation}\label{E:dot_regularized_delta}
    \begin{aligned} 
    \Delta \dot v_\delta &= \dot v_\delta^{(2)} - \dot v_\delta^{(1)} \\
    &= f_{2,\delta}(x_l, x^{(2)}, v_l, v^{(2)}) - f_{1,\delta}(x_l, x^{(1)}, v_l, v^{(1)})
    \end{aligned}
\end{equation}

To analyze this further, we reorganize \eqref{E:dot_regularized_delta} with $f_{1,\delta}(x_l, x^{(2)}, v_l, v^{(2)})$ as follows:
\begin{equation}\label{E:dot_regularized_delta2}
\begin{aligned} 
\Delta \dot v_\delta &= \lb 
\begin{aligned}
f_{2,\delta}(x_l, x^{(2)}, v_l, v^{(2)}) \\
- f_{1,\delta}(x_l, x^{(2)}, v_l, v^{(2)})
\end{aligned}
\rb\\
&\quad + \lb 
\begin{aligned}
f_{1,\delta}(x_l, x^{(2)}, v_l, v^{(2)}) \\
- f_{1,\delta}(x_l, x^{(1)}, v_l, v^{(1)})
\end{aligned}
\rb.
\end{aligned}
\end{equation}

As we investigated previously in \eqref{E:compare_at_same_state}, $f_{2,\delta} \geq f_{1,\delta}$ holds for the entire trajectory of $(x_l, x^{(2)}, v_l, v^{(2)})$ over time $t \in [t_\circ, t_f]$. Therefore, the first term in \eqref{E:dot_regularized_delta2}
\begin{equation} \label{E:compare_at_same_state2}
D \geq \lb \begin{aligned}
f_{2,\delta}(x_l, x^{(2)}, v_l, v^{(2)}) \\
- f_{1,\delta}(x_l, x^{(2)}, v_l, v^{(2)})
\end{aligned}\rb \geq 0.
\end{equation}
where $D$ is a constant that upper-bounds the difference between two models at the same state. 

Now, for the second term in \eqref{E:dot_regularized_delta2}, the term is bounded by the Lipschitz property of $f_{1,\delta}$. That is, there is a constant $L_\delta$ such that 
\begin{equation} \label{E:dot_v_Lipschitz_pre}
\begin{aligned}
    \abs{\Delta \dot v_\delta} = \Bigg | \begin{aligned}
    f_{1,\delta}(x_l, x^{(2)}, v_l, v^{(2)}) \\
    - f_{1,\delta}(x_l, x^{(1)}, v_l, v^{(1)})
    \end{aligned}\Bigg | \leq L_\delta \Big(\abs{\Delta x_\delta} + \abs{\Delta v_\delta}\Big).
    \end{aligned}
\end{equation}

Combining \eqref{E:compare_at_same_state2} and \eqref{E:dot_v_Lipschitz_pre}, we get 
\begin{equation} \label{E:dot_v_Lipschitz_full}
\begin{aligned}
    \abs{\Delta \dot v_\delta} \leq D + L_\delta \Big(\abs{\Delta x_\delta} + \abs{\Delta v_\delta}\Big).
    \end{aligned}
\end{equation}

When we differentiate $\abs{\Delta x_\delta(t)}$ and $\abs{\Delta v_\delta(t)}$, we get
\begin{equation} \label{E:chain_rule_abs}
    \begin{cases}
        \frac{d}{dt}\abs{\Delta x_\delta} \leq \abs{\Delta v_\delta},\\
        \frac{d}{dt}\abs{\Delta v_\delta} \leq \abs{\Delta \dot v_\delta}.
    \end{cases}
\end{equation}

Then, we define $e =\abs{\Delta x_\delta} + \abs{\Delta v_\delta}$, 
\begin{equation} \label{E:chain_rule_abs2}
    \frac{d}{dt} e \leq \abs{\Delta v_\delta} + \abs{\Delta \dot v_\delta}. 
\end{equation}

Based on \eqref{E:dot_v_Lipschitz_full} and \eqref{E:chain_rule_abs2}, we obtain
\begin{equation}\label{E:chain_rule_abs_dot_v_Lipschitz}
\begin{aligned}
\frac{d}{dt} e &\le |\Delta v_\delta| + D + L_\delta \Bigl(|\Delta x_\delta| + |\Delta v_\delta|\Bigr) \\
&\le |\Delta x_\delta| + |\Delta v_\delta| + D + L_\delta  e \\
& = e + L_\delta e + D\\
&= \Bigl(1+L_\delta\Bigr) e + D.
\end{aligned}
\end{equation}

Then, by applying Gronwall's inequality to the inhomogeneous differential inequality, we have 
\begin{equation} \label{E:error_gronwall}
e \leq \frac{D}{1+L_\delta} \Bigl(\exp((1+L_\delta)(t-t_\circ))-1\Bigr)
\end{equation}
with $t \in [t_\circ, t_f]$. Since the initial condition $(x_{l,\circ},x_\circ,v_{l,\circ},v_\circ)$ is the same for $f_{1,\delta}$ and $f_{2,\delta}$ at time $t_\circ$, $e(t_\circ) = 0$. Thus, the structural property in \eqref{E:compare_at_same_state} ensures that initially, the acceleration under $f_{2,\delta}$ is greater than or equal to that under $f_{1,\delta}$. Additionally, \eqref{E:error_gronwall} states that the differences in states between $f_{1, \delta}$ and $f_{2, \delta}$ do not blow up and remain bounded by an exponentially growing function based on Gronwall's inequality. Together, this boundedness and the Lipschitz continuity of the dynamics ensure that the difference maintains its initial sign throughout the finite time interval considered. Consequently, since $f_{2,\delta}$ initially exhibits greater or equal acceleration compared to $f_{1,\delta}$, this inequality persists throughout any finite time interval.

In conclusion, the acceleration $\dot v_\delta^{(2)}$ is always at least as large as $\dot v_\delta^{(1)}$, making $\dot v_\delta^{(2)} - \dot v_\delta^{(1)}$ non-negative for almost every time $t \in [t_\circ, t_f]$ and ultimately $\Delta v_\delta \geq 0$ for almost every time $t \in [t_\circ, t_f]$ since $e(t_\circ) = 0$. 

Given that $\Delta v_\delta \geq 0$ over the time $t \in [t_\circ, t_f]$, we now look at the instantaneous energy cost based on \eqref{E:energy_consumption} which is 
\begin{equation}\label{E:instant_energy_consumption}
    E^{(i)}_\delta= v^{(i)}\,g\bigl(\dot{v}^{(i)}\bigr)
\end{equation}
for $i = 1,2$ with $g$ defined in \eqref{E:g_function}. Since both $\dot v_\delta^{(2)} - \dot v_\delta^{(1)}$ and $v_\delta^{(2)} - v_\delta^{(1)}$ are non-negative for almost every time $t \in [t_\circ, t_f]$, $E^{(1)}_\delta \leq E^{(2)}_\delta$. 

Since energy consumption in \eqref{E:energy_consumption} is defined as the integration of instantaneous energy cost over time interval $[t_\circ, t_f]$ and $E^{(1)}_\delta \leq E^{(2)}_\delta$ for almost every time $t \in [t_\circ, t_f]$, we can claim that energy consumption per unit mass for function $f_{1,\delta}$ and $f_{2,\delta}$, $\omega_{f_{1,\delta}}$ and $\omega_{f_{2,\delta}}$ respectively, have following relationship:
\[
\omega_{f_{1,\delta}} \leq \omega_{f_{2,\delta}}.
\]

Now by letting $\delta\to0$ to make the solutions of regularized function $f_{1,\delta}$ and $f_{2,\delta}$ converge to the classical solution of \eqref{E:dynamics} with $f_{1}$ and $f_{2}$, we get
\[
\lim_{\delta\to0}\omega_{i,\delta}=\omega_i,
\]
for $i=1,2$ where $\omega_i$ is the energy consumption per unit mass for function $f_i$. In conclusion, 
\[
\omega_1 \le \omega_2.
\]

This completes the proof. \end{proof}

\section{Numerical Simulation} \label{sec: numerical_simulation}

\begin{figure}[ht] \vspace{0.2cm}
    \begin{subfigure}[b]{0.45\linewidth}
    \centering
    \includegraphics[width=\linewidth]{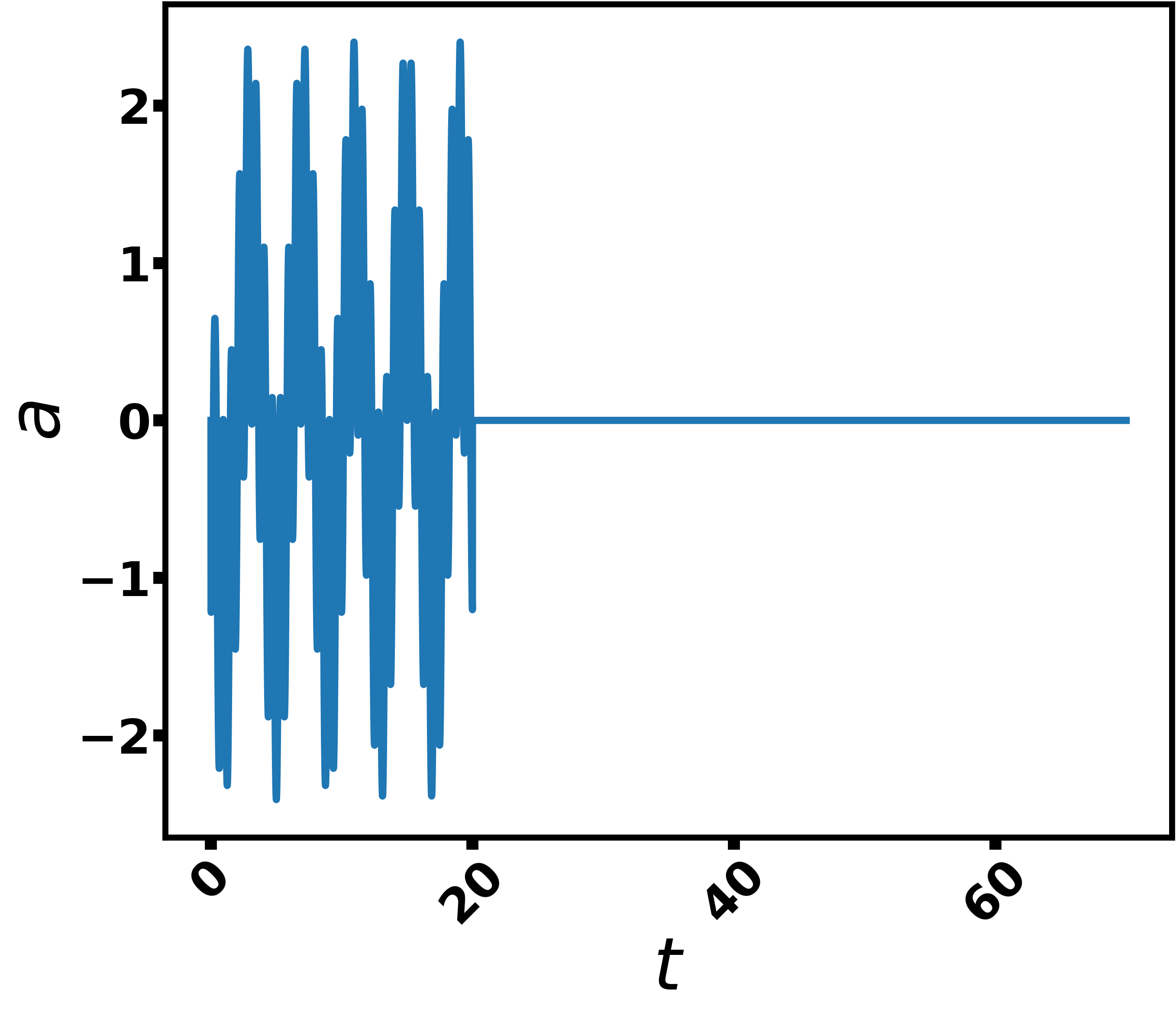} 
    \caption{Acceleration $a_l$} 
    \label{fig: acc4_1} 
  \end{subfigure}
    \begin{subfigure}[b]{0.46\linewidth}
    \centering
    \includegraphics[width=\linewidth]{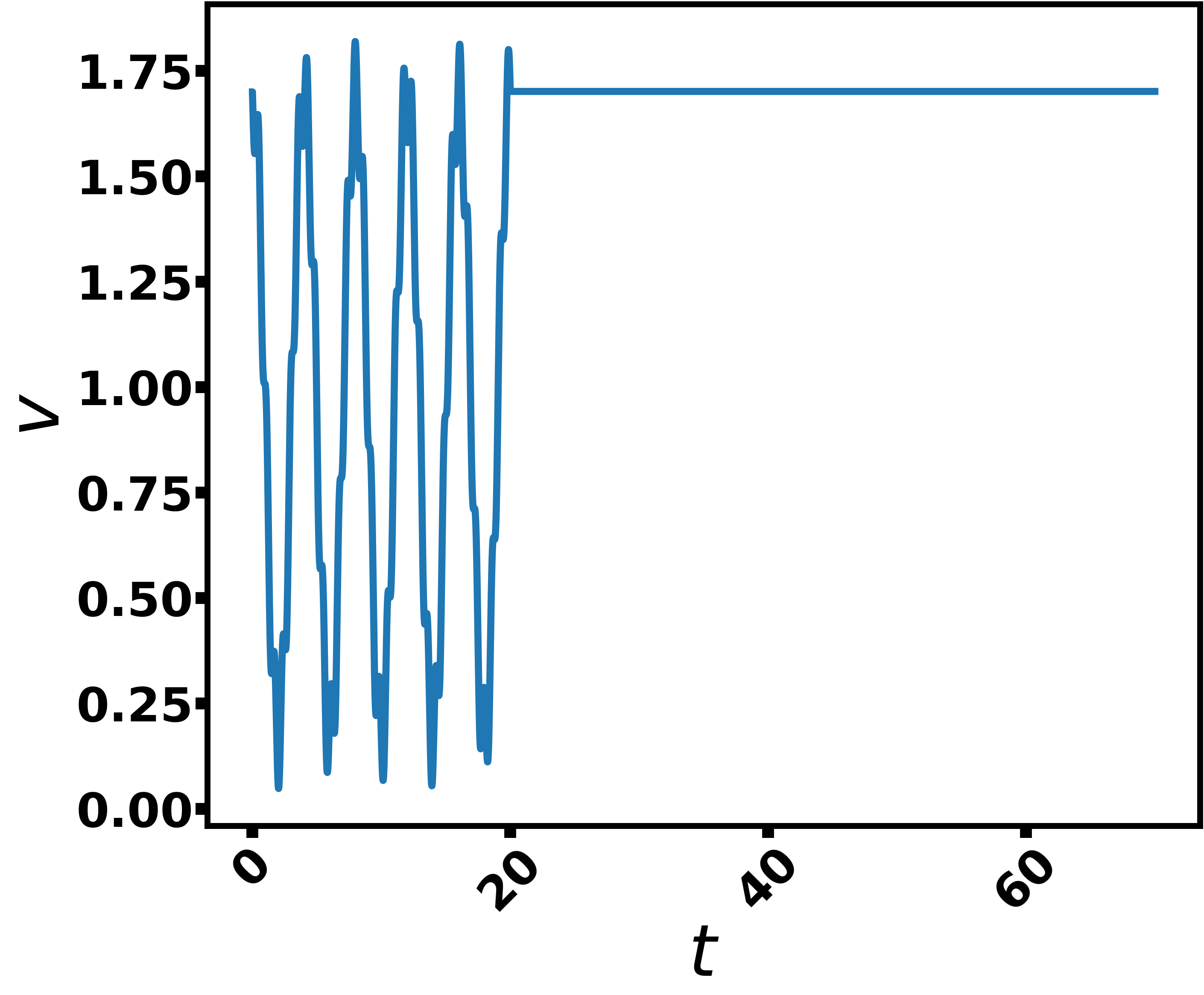} 
    \caption{Velocity $v_l$} 
    \label{fig: vel4_1} 
  \end{subfigure}
  \caption{\small{$a_l$ and $v_l$ in the scenario for the nonlinear stability analysis}}
  \label{fig:acc_vel_4_1} 
\end{figure}

Next, we perform a numerical experiment to evaluate the energy efficiency of the proposed model over a time interval $t \in [0,70]$, under the scenario involving both acceleration and regenerative braking. In this scenario, the leading vehicle, traveling at the velocity of $v_l(t_\circ) = 1.7$ at initial time $t_\circ = 0$, follows a fluctuating acceleration profile defined as 
\begin{equation}\label{E:acc5}
a_l(t) = -\frac{1.65}{1.37} \left[ \sin(0.5\pi t) + \cos(3.2\pi t) \right] \cdot \mathbf{1}_{t \in [0,\,20]}
\end{equation}
where $\mathbf{1}_{t \in [0,\,20]}$ is the indicator function, equal to 1 for $t \in [0,20]$ and 0 otherwise. This acceleration profile generates non-trivial velocity fluctuations, where both energy-consuming and regenerative phases can be observed. The corresponding acceleration $a_l$ and velocity $v_l$ of the leading vehicle in the scenario are shown in Figure \ref{fig:acc_vel_4_1}. 

With the setup of the leading vehicle's acceleration, we consider the simulation of a platoon of 6 vehicles. The leading vehicle ($n = 0$) follows the acceleration defined in \eqref{E:acc5} while the following vehicles $n \in [1,5]$ operate under the proposed model. We denote $x_l$ and $v_l$ for the leading vehicle's position and velocity, while $x_n$ and $v_n$ for the position and velocity of vehicle $n \in [1,5]$. 

At the initial time $t_\circ = 0$, all the following vehicles $n \in [1,5]$ are set to travel at $v_n(t_\circ) = 1.1 v_l(t_\circ)$. The spacing between the leading vehicle and the first following vehicle ($n =1$) is set to $x_l(t_\circ)-x_1(t_\circ) = 0.3$, representing a close-following scenario. The remaining vehicles $n \in [2,5]$ are evenly spaced such that $x_{n-1}(t_\circ)-x_{n}(t_\circ) = 3.5$.

To assess the energy efficiency of the proposed model through the numerical simulation, we run the designed scenario with vehicles in the platoon being operated under the proposed model and under OVFL for $t \in [0, 70]$. For each simulation, we compute the total energy consumption per unit mass of each vehicle $n \in [1, 5]$, with $\eta = 80\%$, over a time interval and perform a comparative analysis. 
\begin{table}[ht!]
    \centering
    \setlength{\tabcolsep}{5pt} 
    \begin{tabular}{|c|c|c|c|c|c|c|} \hline 
         &  n = 1&  n = 2 & n = 3 & n = 4 & n= 5 \\ \thickhline 
         Proposed&  1.827&  1.056&  0.639&  0.487&  0.414\\ \hline 
         OVFL &  1.856&  1.085&  0.662&  0.498&  0.421\\ \hline
         \% Change &  -1.56\% &  -2.67\%&  -3.47\%&  -2.21\%&  -1.66\%\\ \hline
    \end{tabular}
    \caption{\small{Energy consumption per unit mass of Vehicle $n \in [1, 5]$ for the example scenario under Proposed Model and OVFL.}}
    \label{tab:energy_consumption}
\end{table}

Table \ref{tab:energy_consumption} clearly indicates that the total energy consumption per unit mass under the proposed model is consistently lower for each vehicle compared to OVFL. It is important to note that the magnitude of the energy saving depends on the magnitude of $\kappa$. In summary, the numerical simulation supports and validates the analytical result of improved energy efficiency for the proposed model relative to OVFL. 

\section{Stability Analysis} \label{sec:stability}

While the newly introduced control term in \eqref{E: control_term} improves the energy efficiency of the proposed model, the proposed model must exhibit stability under realistic driving conditions to operate safely without unwanted behaviors such as collisions or low velocities. We therefore examine the stability of the proposed model through linear analysis and nonlinear analysis. We also discuss the inherent trade-off between stability and energy efficiency of the proposed model due to the newly introduced control term in \eqref{E: control_term}.

We first examine the linear local stability of the proposed model near the equilibrium. We reformulate \eqref{E:main_dynamics} into a relative framework that captures the deviation between the ego vehicle and the leading vehicle. For assessing the local stability of the proposed model, we focus on a steady-state scenario where a leading vehicle is moving at a constant speed ($v_l = \bar v$ and $a_l = 0$). Then, the resulting system of the spacing ($z = x_l - x$) and the relative velocity ($y = v_l - v$) at the steady-state scenario is
\begin{equation*}\label{E:main_dynamics_difference_eq}
\begin{cases}
\dot z  = y, \\
\dot y = - \left( \alpha \lb V(z)-v_l+y \rb + \beta\frac{y}{z^2} -  \kappa \frac{(v_l-y)^2 y^2}{y^2+\epsilon} \right),  \\
(z(t_\circ) , y(t_\circ)) = (z_\circ, y_\circ)
\end{cases}
\end{equation*}
where the equilibrium spacing ($z_{eq}$) and relative velocity ($y_{eq}$) are
\begin{equation}\label{E:equilibrium}
(z_{eq},y_{eq}) = (V^{-1}(\bar v),0).
\end{equation}

The proposed model has the same equilibrium as the OVM, and its linearized dynamics around the equilibrium are also the same. As a result, the linear local stability characteristics are identical and have been extensively studied in \cite{wilson2011car, ward2011criteria}.

With the linear stability analysis established, we proceed to evaluate the nonlinear stability of the proposed model through numerical simulation performed in Section \ref{sec: numerical_simulation}. The simulation setup enables us to examine how vehicles in the proposed model respond to the leading vehicle's fluctuating motion. Specifically, we evaluate the model’s ability to suppress perturbation growth, maintain stability, converge to equilibrium states, and prevent undesired outcomes such as collisions or negative velocities.

\begin{figure}[ht!] 
    \centering
    \begin{subfigure}[b]{0.49\linewidth}
    \centering
    \includegraphics[width=\linewidth]{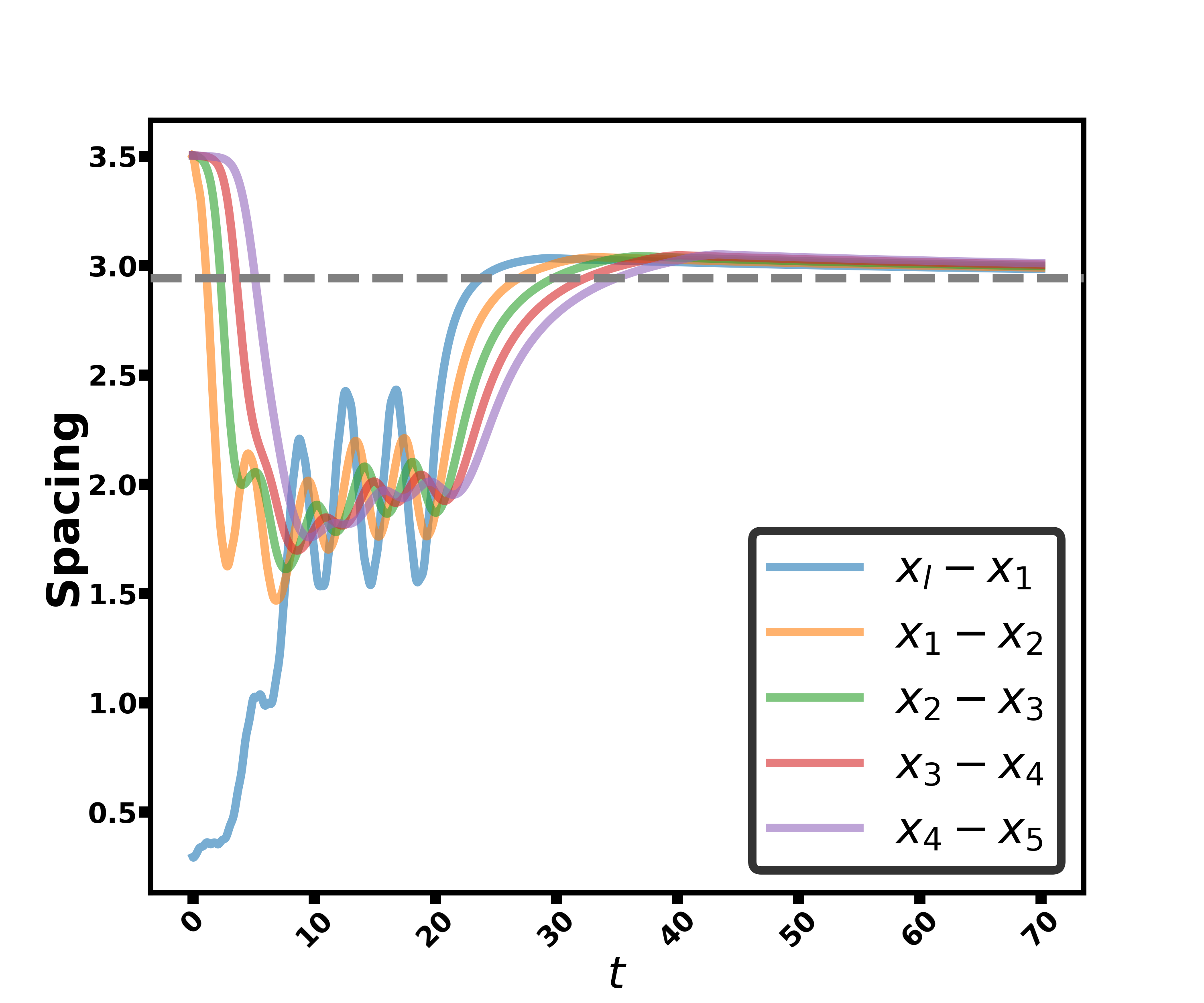} 
    \caption{Spacing} 
  \end{subfigure}
    \begin{subfigure}[b]{0.49\linewidth}
    \centering
    \includegraphics[width=\linewidth]{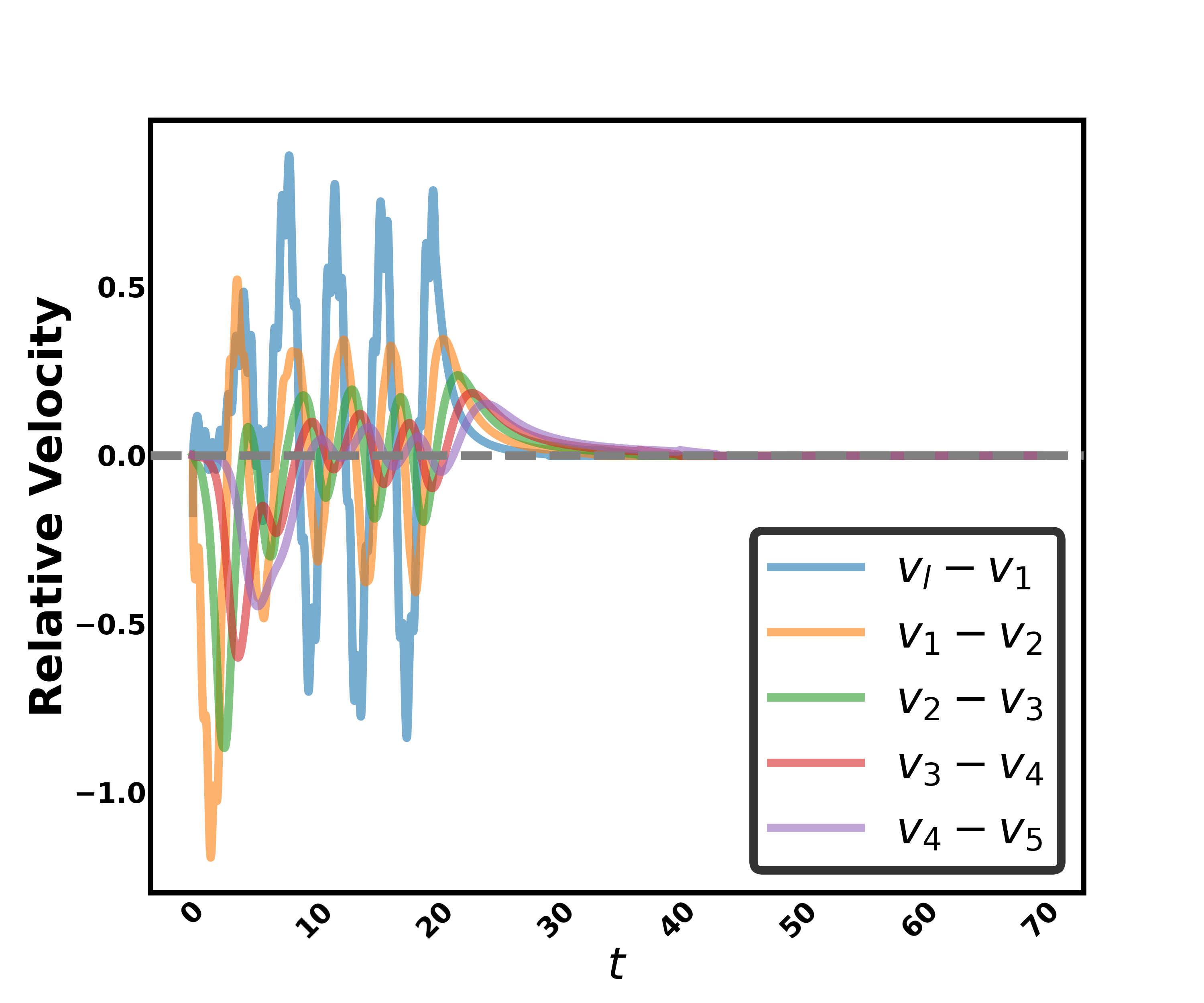} 
    \caption{Relative Velocity} 
  \end{subfigure}
  \caption{\small{Spacing and Relative Velocity between vehicles in the platoon under the proposed model for the numerical simulation in Section \ref{sec: numerical_simulation}. Gray Dashed lines are the equilibrium spacing and relative velocity based on $v_l$ at $t = 20$.}}
  \label{fig:proposed_nonlinear} 
\end{figure}

Figure \ref{fig:proposed_nonlinear} presents the time evolution of spacing and relative velocity between vehicles over $t \in [0, 70]$.  

The simulation results reveal that the magnitude of perturbation due to the leading vehicle's motion attenuates as it propagates through the platoon. Over time, the follower vehicles in the platoon converge to the expected equilibrium spacing and relative velocity. Therefore, the numerical experiment corroborates the linear local stability proof. Moreover, throughout the simulation, no instance of collision or negative velocity is observed. 

Building on the linear stability analysis and supporting simulation, we now discuss the trade-off between the energy efficiency and stability of the proposed model due to the newly introduced control term stated in \eqref{E: control_term}. From \eqref{E:main_dynamics}, we know that the control term encourages deceleration, thereby lowering the velocity of the vehicle. While this contributes to improving the energy efficiency of the vehicle, it slows down the rate of convergence to the equilibrium.

In addition, the control term destabilizes the system. Either when $v$ is sufficiently large or the magnitude of $\kappa$ is large, then the control term in \eqref{E: control_term} can dominate the other two terms in \eqref{E:main_dynamics}. In such cases, the velocity goes to zero, bringing the vehicles to a stop. This indicates a loss of global stability with respect to the non-zero speed equilibrium in \eqref{E:equilibrium}.  

Therefore, the value of $\kappa$ must be chosen small enough to ensure a reasonable rate of convergence and to preserve the global stability of the model, while also minimizing the magnitude of $I$ to improve the energy efficiency.

\section{Conclusions/Future Works}

Considering battery dynamics under reasonable assumptions, this paper introduced an energy-efficient nonlinear dynamical model for A-EVs. The improved energy efficiency of the proposed model is demonstrated both analytically and numerically. The stability properties are examined, along with the trade-off between energy efficiency and stability.

The potential future work includes 1) relaxing the assumptions made in the current model to more accurately reflect the battery dynamics, 2) performing analytical investigations of nonlinear stability, collision avoidance, and negative velocity behaviors, and 3) determining the optimal value of $\kappa$ to effectively balance the energy efficiency and the stability.

\bibliographystyle{plain} 
\bibliography{reference.bib} 
\end{document}